# Eikonal slant helices and eikonal Darboux helices in 3-dimensional pseudo-Riemannian manifolds

**Mehmet Önder[a], Evren Zıplar[b]**

[a]*Celal Bayar University, Faculty of Arts and Sciences, Department of Mathematics, Muradiye Campus, 45047 Muradiye, Manisa, Turkey.*
E-mail: mehmet.onder@cbu.edu.tr

[b]*Çankırı Karatekin University, Faculty of Science, Department of Mathematics, Çankırı, Turkey*
E-mail: evrenziplar@karatekin.edu.tr

**Abstract**

In this study, we give definitions and characterizations of eikonal slant helices, eikonal Darboux helices and non-normed eikonal Darboux helices in 3-dimensional pseudo-Riemannian manifold $M$. We show that every eikonal slant helix is also an eikonal Darboux helix for timelike and spacelike curves. Furthermore, we obtain that if the non-null curve $\alpha$ is a non-normed eikonal Darboux helix, then $\alpha$ is an eikonal slant helix if and only if $\varepsilon_3 \kappa^2 + \varepsilon_1 \tau^2 = \text{constant}$, where $\kappa$ and $\tau$ are curvature and torsion of $\alpha$, respectively. Finally, we define null-eikonal helices, slant helices and Darboux helices. Also, we give their characterizations.

## 1. Introduction

In the nature and science, some special curves have an important role and many applications. The well-known of such curves is helix curve. In the Euclidean 3-space $E^3$, a general helix is defined as a special curve whose tangent line makes a constant angle with a fixed straight line which is called the axis of the helix [4]. This definition gives that the tangent indicatrix of a general helix is a planar curve. Moreover, the classical result for the helices first was given by Lancret in 1802 and proved by B. de Saint Venant in 1845 as follows: *A necessary and sufficient condition that a curve to be a general helix is that the ratio of the first curvature to the second curvature be constant i.e., $\kappa/\tau$ is constant along the curve, where $\kappa$ and $\tau$ denote the first and second curvatures of the curve, respectively* [20]. The same definition is also valid in Lorentzian space and spacelike, timelike and null helices have been studied by some mathematicians [7-9].

Furthermore, there exist more special curves in the space such as slant helix which first introduced by Izumiya and Takeuchi by the property that the normal lines of curve make a constant angle with a fixed direction in the Euclidean 3-space $E^3$ [14]. Slant helices have been studied by some mathematicians and new kinds of these curves also have been introduced [1,11,16,17,19]. Moreover, these curves have been considered in Lorentzian spaces [2,3].

Later, a new kind of helices has been defined by Zıplar, Şenol and Yaylı according to the Darboux vector of a space curve in $E^3$. They have called this new curve as Darboux helix which is defined by the property that the Darboux vector of a space curve makes a constant angle with a fixed direction and they have given the characterizations of this new special curve [22].

Let $M$ be a Riemannian manifold with the metric $g$ and $f: M \to \mathbb{R}$ be a function with gradient $\nabla f$. The function $f$ is called eikonal if $\|\nabla f\|$ is constant [5]. There exist many applications of $\nabla f$ in mathematical physics and geometry. For instance, if $f$ is non-constant on connected $M$, then the Riemannian condition $\|\nabla f\|^2 = 1$ is precisely the eikonal equation of geometrical optics. So, on a connected $M$, a non-constant real valued function $f$ is Riemannian if $f$ satisfies this eikonal equation. In the geometrical optical interpretation, the level sets of $f$ are interpreted as wave fronts. The characteristics of the eikonal equation (as a partial differential equation), are then the solutions of the gradient flow equation for $f$ (an ordinary differential equation), $x' = \nabla f$, which are geodesics of $M$ orthogonal to the level sets of $f$, and which are parameterized by arc length. These geodesics can be interpreted as light rays orthogonal to the wave fronts (See [10] for details). Later, Şenol, Zıplar and Yaylı have defined eikonal helices and eikonal slant helices by considering a space curve with a function $f: M \to \mathbb{R}$ where $M$ is a Riemannian manifold [21].

In this study, we define and give the characterizations of $f$-eikonal slant helices and $f$-eikonal Darboux helices for non-null and null curves in a pseudo-Riemannian manifold. For this purpose, we need the following definitions.

**Definition 1.1.** ([18]) A metric tensor $g$ in a smooth manifold $M$ is a symmetric non-degenerate $(0,2)$ tensor field in $M$.

On the other hand if $TM$ is the tangent bundle of $M$, then for all $X, Y \in TM$, $g(X,Y) = g(Y,X)$ and at each point $p$ of $M$, if $g(X_p, Y_p) = 0$ for all $Y_p \in T_p(M)$, then $X_p = 0$ (non-degenerate) where $T_p(M)$ is the tangent space of $M$ at the point $p$ and $g: T_p(M) \times T_p(M) \to \mathbb{R}$.

**Definition 1.2.** ([18]) A pseudo-Riemannian manifold (or semi-Riemannian manifold) is a smooth manifold $M$ furnished with a metric tensor $g$. That is, a pseudo-Riemannian manifold is an ordered pair $(M, g)$.

**Definition 1.3.** Let $M$ be a pseudo-Riemannian manifold and $g$ be its metric. For the function $f: M \to \mathbb{R}$, it is said that $f$ is eikonal if $\|\nabla f\|$ is constant, where $\nabla f$ is gradient of $f$, i.e., $df(X) = g(\nabla f, X) = X(f)$.

**Lemma 1.1.** ([18]) Let $(M,g)$ be a pseudo-Riemannian manifold and $\nabla$ be the Levi-Civita connection of $M$. The Hessian $H^f$ of a $f \in F(M)$ is the symmetric (0,2) tensor field such that

$$H^f(X,Y) = g(\nabla_X (\text{grad} f), Y),$$

where $F(M)$ shows the set of differentiable functions defined on $M$.

From Lemma 1.1, we have the following corollary.

**Corollary 1.1.** The Hessian $H^f$ of a $f \in F(M)$ is zero, i.e., $H^f = 0$ if and only if $\nabla f$ is parallel in $M$.

## 2. Non-null Eikonal Slant Helices and Non-null Eikonal Darboux Helices

Let $(M,g)$ be a time-oriented 3-dimensional pseudo-Riemannian manifold and $\alpha : I \to M$ be a unit speed curve on $M$, i.e., $g(\alpha', \alpha') = \varepsilon_1 = \pm 1$ is satisfied along $\alpha$ where $\alpha'$ is the velocity vector filed of the curves and $g$ shows the metric tensor (or Lorentzian metric) given by $g(a,b) = -a_1 b_1 + a_2 b_2 + a_3 b_3$ for the vectors $a = (a_1, a_2, a_3)$, $b = (b_1, b_2, b_3) \in TM$. The constant $\varepsilon_1 = \pm 1$ defined by $\varepsilon_1 = g(\alpha', \alpha')$ is called the causal character of $\alpha$. Then, a unit speed curve $\alpha$ is said to be spacelike or timelike if its causal character is 1 or -1, respectively. The curve $\alpha$ is said to be a Frenet curve if $g(\alpha'', \alpha'') \neq 0$. Like Euclidean geometry, every Frenet curve $\alpha$ on $(M,g)$ admits an orthonormal Frenet frame field $\{V_1, V_2, V_3\}$ along $\alpha$ such that $V_1 = \alpha'(s)$. The vector fields $V_1, V_2, V_3$ are called tangent vector field, principle normal vector field and binormal vector field of $\alpha$, respectively and $\{V_1, V_2, V_3\}$ satisfies the following Frenet-Serret formula:

$$\begin{bmatrix} \nabla_{V_1} V_1 \\ \nabla_{V_1} V_2 \\ \nabla_{V_1} V_3 \end{bmatrix} = \begin{bmatrix} 0 & \varepsilon_2 \kappa & 0 \\ -\varepsilon_1 \kappa & 0 & -\varepsilon_3 \tau \\ 0 & \varepsilon_2 \tau & 0 \end{bmatrix} \begin{bmatrix} V_1 \\ V_2 \\ V_3 \end{bmatrix},$$

where $\nabla$ is the Levi-Civita connection of $(M,g)$ [12,13,15]. The functions $\kappa \geq 0$ and $\tau$ are called the curvature and torsion, respectively. The constants $\varepsilon_2$ and $\varepsilon_3$ are defined by

$$\varepsilon_i = g(V_i, V_i), \ i = 2,3.$$

and called second causal character and third causal character of $\alpha$, respectively. Note that $\varepsilon_3 = -\varepsilon_1 \varepsilon_2$ and $V_i \times V_j = \varepsilon_i \varepsilon_j V_k$, where $(i,j,k) = (1,2,3), (2,3,1), (3,1,2)$.

The vector $W = \tau V_1 - \kappa V_3$ is called Darboux vector of the curve $\alpha$. Then for the Frenet formulae we have $\nabla_{V_1} V_i = W \times V_i$, $(i = 1,2,3)$; where "$\times$" shows the vector product in $M$.

As in the case of Riemannian geometry, a Frenet curve $\alpha$ is a geodesic if and only if $\kappa = 0$. A circular helix is a Frenet curve whose curvature and torsion are constants. If the curvature $\kappa$ is constant and the torsion $\tau$ is zero, then the curve is called a pseudo circle.

Pseudo circles are regarded as degenerate helices. Helices, which are not circles, are frequently called proper helices.

**Definition 2.1.** Let $M^3$ be a 3-dimensional pseudo-Riemannian manifold with the Lorentzian metric $g$ and let $\alpha(s)$ be a non-null Frenet curve with the Frenet frame $\{V_1, V_2, V_3\}$ in $M^3$. Let $f: M^3 \to \mathbb{R}$ be an eikonal function along curve $\alpha$, i.e. $\|\nabla f\| = $ constant along the curve $\alpha$. If the function $g(\nabla f, V_2)$ is a non-zero constant along $\alpha$, then $\alpha$ is called a non-null $f$-eikonal slant helix. And, $\nabla f$ is called the axis of the $f$-eikonal slant helix $\alpha$.

**Definition 2.2.** Let $M^3$ be a pseudo-Riemannian manifold with the Lorentzian metric $g$ and $\alpha$ be a non-null Frenet curve in $M^3$ with Frenet frame $\{V_1, V_2, V_3\}$, non-zero curvatures $\kappa, \tau$ and Darboux vector $W = \tau V_1 - \kappa V_3$. Also, let $f: M^3 \to \mathbb{R}$ be an eikonal function along $\alpha$. If the unit Darboux vector

$$W_0 = \frac{\tau}{\sqrt{|\varepsilon_3 \kappa^2 + \varepsilon_1 \tau^2|}} V_1 - \frac{\kappa}{\sqrt{|\varepsilon_3 \kappa^2 + \varepsilon_1 \tau^2|}} V_3,$$

of the curve $\alpha$ makes a constant angle $\varphi$ with the gradient of the function $f$, that is $g(W_0, \nabla f)$ is constant along $\alpha$, then the curve $\alpha$ is called a non-null $f$-eikonal Darboux helix.

Especially, if $g(W, \nabla f) = $ constant, then $\alpha$ is called a non-normed non-null $f$-eikonal Darboux helix. Then, we have the following Corollary.

***Corollary 2.1.*** *A non-normed non-null $f$-eikonal Darboux helix is a non-null $f$-eikonal Darboux helix if and only if $\varepsilon_3 \kappa^2 + \varepsilon_1 \tau^2$ is constant.*

**Example 2.1.** We consider the pseudo-Riemannnian manifold $M^3 = \mathbb{R}_1^3$ with the Lorentzian metric $g$. Let

$$f: M^3 \to \mathbb{R}$$
$$(x, y, z) \to f(x, y, z) = x^2 + y^2 + z$$

be a function defined in $M^3$ and consider the spacelike curve

$$\alpha : I \subset \mathbb{R} \to M^3$$

$$s \to \alpha(s) = \left( a \cosh \frac{s}{\sqrt{a^2 + b^2}}, \, a \sinh \frac{s}{\sqrt{a^2 + b^2}}, \, \frac{bs}{\sqrt{a^2 + b^2}} \right); \quad a, b > 0$$

in $M^3$. If we compute $\nabla f$, we find out $\nabla f$ as $\nabla f = (2x, 2y, 1)$. Then, we have

$$\|\nabla f\| = \sqrt{\left|1 + 4\left(-x^2 + y^2\right)\right|},$$

and, along the curve $\alpha$, we find out

$$\|\nabla f\| = \sqrt{\left|1 - 4a^2\right|} = \text{constant}.$$

That is, $f$ is an eikonal function along $\alpha$. Moreover, by a simple computation we have that the principal normal of the curve is

$$V_2(s) = \left( \cosh \frac{s}{\sqrt{a^2 + b^2}}, \, \sinh \frac{s}{\sqrt{a^2 + b^2}}, \, 0 \right).$$

Since

$$\nabla f = \left( 2a \cosh \frac{s}{\sqrt{a^2 + b^2}}, \, 2a \sinh \frac{s}{\sqrt{a^2 + b^2}}, \, 1 \right),$$

along $\alpha$, we easily see that $g(\nabla f, V_2) = -2a = \text{constant}$ which means that $\alpha$ is a non-null $f$-eikonal slant helix in $M^3$.

On the other hand, non-normed Darboux vector of $\alpha$ is

$$W = \left( -\frac{a(a^2 + b^2)^2 + ab^2}{b(a^2 + b^2)^{3/2}} \sinh \frac{s}{\sqrt{a^2 + b^2}}, \, -\frac{a(a^2 + b^2)^2 + ab^2}{b(a^2 + b^2)^{3/2}} \cosh \frac{s}{\sqrt{a^2 + b^2}}, \right.$$
$$\left. \frac{a^2 - (a^2 + b^2)^2}{(a^2 + b^2)^{3/2}} \right)$$

and curvatures are $\kappa = \frac{a}{a^2 + b^2}$, $\tau = -\frac{a^2 + b^2}{b}$, respectively. Then we obtain that

$$g(\nabla f, W) = \frac{a^2 - (a^2 + b^2)^2}{(a^2 + b^2)^{3/2}} = \text{constant},$$

along $\alpha$. So, $\alpha$ is a non-null $f$-eikonal non-normed Darboux helix curve in $M^3$. Since $\kappa, \tau$ are constants $\alpha$ is also a non-null $f$-eikonal Darboux helix curve in $M^3$.

Now, we give some theorems concerned with non-null $f$-eikonal slant helices and $f$-eikonal Darboux helices in pseudo-Riemannian manifold. Whenever we write $M^3$, we will consider $M^3$ as a 3-dimensional pseudo-Riemannian manifold with the Lorentzian metric $g$.

**Theorem 2.1.** *Let $\alpha : I \subset \mathbb{R} \to M^3$ be a non-null curve in $M^3$ with non-zero curvatures $\kappa, \tau$ and assume that $\alpha(s)$ is not a helix. Let $f : M^3 \to \mathbb{R}$ be an eikonal function along curve $\alpha$ and the Hessian $H^f = 0$. If $\alpha(s)$ is a non-null f-eikonal slant helix curve in $M^3$, then the following properties hold:*

*i) The function*

$$\frac{\kappa^2}{\left(\varepsilon_1 \tau^2 + \varepsilon_3 \kappa^2\right)^{3/2}} \left(\frac{\tau}{\kappa}\right)', \tag{1}$$

*is a real constant.*

*ii) The axis of $f$-eikonal slant helix is obtained as*

$$\nabla f = \frac{n\tau}{\sqrt{\left|\varepsilon_1 \tau^2 + \varepsilon_3 \kappa^2\right|}} V_1 + c V_2 - \frac{n\kappa}{\sqrt{\left|\varepsilon_1 \tau^2 + \varepsilon_3 \kappa^2\right|}} V_3,$$

*where $c$ and $n$ are non-zero constants.*

**Proof. i)** Since $\alpha$ is a non-null $f$-eikonal slant helix, we have $g(\nabla f, V_2) = c = $ constant. So, there exist smooth functions $a_1 = a_1(s)$, $a_2 = a_2(s) = c$ and $a_3 = a_3(s)$ of arc length $s$ such that

$$\nabla f = a_1 V_1 + c V_2 + a_3 V_3, \tag{2}$$

where $\{V_1, V_2, V_3\}$ is a basis of $TM^3$ (tangent bundle of $M^3$).

From Corollary 1.1, $\nabla f$ is parallel in $M^3$, i.e., $\nabla_{V_1} \nabla f = 0$ along $\alpha$. Then, if we take the derivative in each part of (2) in the direction $V_1$ in $M^3$ and use the Frenet equations, we get

$$(V_1[a_1] - \varepsilon_1 \kappa c) V_1 + (\varepsilon_2 a_1 \kappa + \varepsilon_2 a_3 \tau) V_2 + (V_1[a_3] - \varepsilon_3 \tau c) V_3 = 0, \tag{3}$$

Since $V_1[a_i] = a_i'(s)$, $(i = 1, 2, 3)$ in (3) and the Frenet frame $\{V_1, V_2, V_3\}$ is linearly independent, we have

$$\begin{cases} a_1' - \varepsilon_1 \kappa c = 0, \\ a_1 \kappa + a_3 \tau = 0, \\ a_3' - \varepsilon_3 \tau c = 0. \end{cases} \tag{4}$$

From the second equation of the system (4) we obtain

$$a_1 = -\left(\frac{\tau}{\kappa}\right)a_3. \tag{5}$$

Since $f$ is an eikonal function along $\alpha$, we have $\|\nabla f\|$ is constant. Then (2) and (5) give that

$$\left[\varepsilon_1\left(\frac{\tau}{\kappa}\right)^2 + \varepsilon_3\right]a_3^2 + \varepsilon_2 c^2 = \text{constant}, \tag{6}$$

and from (6) we can write

$$\left[\varepsilon_1\left(\frac{\tau}{\kappa}\right)^2 + \varepsilon_3\right]a_3^2 = n^2, \tag{7}$$

where $n^2$ is a constant. Since $\alpha$ is not a helix curve in $M^3$ and curvatures are not zero, we have that $n$ is a non-zero constant. Then, from (7) we have

$$a_3 = \pm \frac{n}{\sqrt{\left|\varepsilon_1\left(\frac{\tau}{\kappa}\right)^2 + \varepsilon_3\right|}}. \tag{8}$$

By taking the derivative of (8) with respect to $s$ and using the third equation of the system (4), we get that the function

$$\frac{\kappa^2}{\left(\varepsilon_1\tau^2 + \varepsilon_3\kappa^2\right)^{3/2}}\left(\frac{\tau}{\kappa}\right)', \tag{9}$$

is a constant, which is desired function.

*ii*) By direct calculation from (5) and (8), we have

$$a_1 = \frac{n\tau}{\sqrt{\left|\varepsilon_1\tau^2 + \varepsilon_3\kappa^2\right|}}, \quad a_3 = \frac{n\kappa}{\sqrt{\left|\varepsilon_1\tau^2 + \varepsilon_3\kappa^2\right|}},$$

where $n$ is a non-zero constant. Then, from (2) the axis of $f$-eikonal slant helix is

$$\nabla f = \frac{n\tau}{\sqrt{\left|\varepsilon_1\tau^2 + \varepsilon_3\kappa^2\right|}}V_1 + cV_2 - \frac{n\kappa}{\sqrt{\left|\varepsilon_1\tau^2 + \varepsilon_3\kappa^2\right|}}V_3. \tag{10}$$

The above Theorem has the following corollary.

***Corollary 2.2.*** *Let $\alpha: I \subset \mathbb{R} \to M^3$ be a non-null curve in $M^3$ with non-zero curvatures $\kappa, \tau$ and assume that $\alpha(s)$ is not a helix. Let $f: M^3 \to \mathbb{R}$ be an eikonal function along curve $\alpha$* and the Hessian $H^f = 0$. If $\alpha(s)$ is a non-null f-eikonal slant helix curve in $M^3$, then, the curvatures $\kappa$ and $\tau$ satisfy the following non-linear equation system:

$$\left(\frac{n\tau}{\sqrt{|\varepsilon_1\tau^2 + \varepsilon_3\kappa^2|}}\right)' - \varepsilon_1\kappa c = 0, \quad \left(\frac{n\kappa}{\sqrt{|\varepsilon_1\tau^2 + \varepsilon_3\kappa^2|}}\right)' - \varepsilon_3\tau c = 0. \tag{11}$$

**Theorem 2.2.** *Let $\alpha : I \subset \mathbb{R} \to M^3$ be a non-null curve in $M^3$ with non-zero curvatures $\kappa, \tau$ and assume that $\alpha(s)$ is not a helix. Let $f : M^3 \to \mathbb{R}$ be an eikonal function along curve $\alpha$ and the Hessian $H^f = 0$. Then, every non-null f-eikonal slant helix in $M^3$ is also a non-null $f$-eikonal Darboux helix in $M^3$.*

**Proof.** Let $\alpha$ be a non-null $f$-eikonal slant helix in $M^3$. Then, from Theorem 2.1, the axis of $\alpha$ is

$$\nabla f = \frac{n\tau}{\sqrt{|\varepsilon_1\tau^2 + \varepsilon_3\kappa^2|}} V_1 + cV_2 - \frac{n\kappa}{\sqrt{|\varepsilon_1\tau^2 + \varepsilon_3\kappa^2|}} V_3. \tag{12}$$

Considering the unit Darboux vector $W_0$, equality (12) can be written as follows

$$\nabla f = nW_0 + cV_2, \tag{13}$$

which shows that $\nabla f$ lies on the plane spanned by $W_0$ and $V_2$. Since $n$ is a non-zero constant, from (13), we have $g(\nabla f, W_0) = n$ is constant along $\alpha$, i.e, $\alpha$ is a non-null $f$-eikonal Darboux helix in $M^3$.

**Theorem 2.3.** *Let $\alpha : I \subset \mathbb{R} \to M^3$ be a non-null curve in $M^3$ with non-zero curvatures $\kappa, \tau$ and assume that $\alpha(s)$ is not a helix. Let $f : M^3 \to \mathbb{R}$ be an eikonal function along curve $\alpha$ and the Hessian $H^f = 0$. Let $\alpha$ be a non-normed non-null $f$-eikonal Darboux helix with Darboux vector $W$. Then $\alpha$ is a non-null $f$-eikonal slant helix if and only if $\|W\|$ is a non-zero constant.*

**Proof.** Since $\alpha$ is a non-normed non-null $f$-eikonal Darboux helix, we have $g(W, \nabla f) = \text{constant}$. On the other hand, there exist smooth functions $a_1 = a_1(s)$, $a_2 = a_2(s)$ and $a_3 = a_3(s)$ of arc length $s$ such that

$$\nabla f = a_1 V_1 + a_2 V_2 + a_3 V_3, \tag{14}$$

where $a_1, a_2, a_3$ are assumed non-zero and $\{V_1, V_2, V_3\}$ is a basis of $TM^3$. From Corollary 1.1., $\nabla f$ is parallel in $M^3$, i.e., $\nabla_{V_1}\nabla f = 0$ along $\alpha$. Then, if we take the derivative in each part of (14) in the direction $V_1$ in $M^3$ and use the Frenet equations, we get

$$(a_1' - \varepsilon_1 a_2 \kappa)V_1 + (\varepsilon_2 a_1 \kappa + a_2' + \varepsilon_2 a_3 \tau)V_2 + (a_3' - \varepsilon_3 a_2 \tau)V_3 = 0, \tag{15}$$

where $a_i'(s) = V_1[a_i]$, $(i = 1, 2, 3)$. Since the Frenet frame $\{V_1, V_2, V_3\}$ is linearly independent, we have

$$\begin{cases} a_1' - \varepsilon_1 \kappa a_2 = 0, \\ a_2' + \varepsilon_2 \kappa a_1 + \varepsilon_2 \tau a_3 = 0, \\ a_3' - \varepsilon_3 \tau a_2 = 0. \end{cases} \tag{16}$$

Equality $g(W, \nabla f) = \text{constant}$ gives that

$$\varepsilon_1 a_1 \tau - \varepsilon_3 a_3 \kappa = \text{constant}. \tag{17}$$

Differentiating (17) and using the first and third equations of system (16) we obtain

$$\varepsilon_1 a_1 \tau' - \varepsilon_3 a_3 \kappa' = 0 \tag{18}$$

From (18) and the second equation of system (16) it follows

$$a_2' = -\frac{\varepsilon_2 (\varepsilon_1 \varepsilon_3 \kappa^2 + \tau^2)'}{2\tau'} a_3 \tag{19}$$

In (19) if $a_3 = 0$, from (16) we have $a_1 = a_2 = 0$, i.e., $\nabla f = 0$ which is a contradiction. Then we have $a_3 \neq 0$ and from (19) we see that $a_2 = a_2(s)$ is constant if and only if $\varepsilon_1 \varepsilon_3 \kappa^2 + \tau^2 = \text{constant}$ which means that $\varepsilon_1 \tau^2 + \varepsilon_3 \kappa^2 = \text{constant}$, i.e, $\alpha$ is a non-null $f$-eikonal slant helix if and only if $\|W\|$ is a non-zero constant.

From Theorem 2.3 and Corollary 2.1, we have the following corollary.

***Corollary 2.3.*** *Let $\alpha : I \subset \mathbb{R} \to M^3$ be a non-null curve in $M^3$ with non-zero curvatures $\kappa, \tau$ and assume that $\alpha(s)$ is not a helix. Let $f : M^3 \to \mathbb{R}$ be an eikonal function along curve $\alpha$ and the Hessian $H^f = 0$. Let $\alpha$ be a non-normed non-null $f$-eikonal Darboux helix. Then $\alpha$ is a non-null $f$-eikonal slant helix if and only if $\alpha$ is a non-null $f$-eikonal Darboux helix.*

## 3. Null Eikonal Helices, Slant Helices and Darboux Helices

Let $\alpha$ be a curve in 3-dimensional pseudo-Riemannian manifold $(M, g)$. Then, the curve $\alpha$ is called a null curve if $g(V_1, V_1) = 0$. By a Cartan frame or Frenet frame $\{V_1, V_2, V_3\}$ of $\alpha$, we mean a family of vector fields $V_1 = V_1(s)$, $V_2 = V_2(s)$, $V_3 = V_3(s)$ along the curve $\alpha$ satisfying the following conditions:

$$\begin{cases} \alpha'(s) = V_1, \ g(V_1,V_1) = g(V_2,V_2) = 0, \ g(V_1,V_2) = 1, \\ g(V_1,V_3) = g(V_2,V_3) = 0, \ g(V_3,V_3) = 1, \\ V_1 \times V_2 = V_3, \ V_2 \times V_3 = V_2, \ V_3 \times V_1 = V_1. \end{cases} \quad (20)$$

([6]). Here $V_1$, $V_2$ and $V_3$ are called tangent vector field, binormal vector field and (principal) normal vector field of $\alpha$, respectively. Then the derivative formula of the frame is given as follows

$$\begin{bmatrix} \nabla_{V_1} V_1 \\ \nabla_{V_1} V_2 \\ \nabla_{V_1} V_3 \end{bmatrix} = \begin{bmatrix} 0 & 0 & \kappa \\ 0 & 0 & \tau \\ -\tau & -\kappa & 0 \end{bmatrix} \begin{bmatrix} V_1 \\ V_2 \\ V_3 \end{bmatrix}, \quad (21)$$

where $\kappa$ and $\tau$ are called the curvature and torsion of $\gamma$, respectively [6].

The vector $W = \tau V_1 - \kappa V_2$ is called Darboux vector of the curve $\alpha$. Then for the Frenet formulae (21) we have $\nabla_{V_1} V_i = W \times V_i$, $(i = 1, 2, 3)$; where "$\times$" shows the vector product in $M^3$.

**Definition 3.1.** Let $M^3$ be a 3-dimensional pseudo-Riemannian manifold with the metric $g$ and let $\alpha(s)$ be a null Frenet curve with the Frenet frame $\{V_1, V_2, V_3\}$ and Darboux vector $W = \tau V_1 - \kappa V_2$ in $M^3$. Let $f : M^3 \to \mathbb{R}$ be an eikonal function along the curve $\alpha$, i.e. $\|\nabla f\| = $ constant along $\alpha$. Then we define the followings,

**i)** If the function $g(\nabla f, V_1)$ is a non-zero constant along $\alpha$, then $\alpha$ is called a null $f$-eikonal helix curve. And, $\nabla f$ is called the axis of the null $f$-eikonal helix curve $\alpha$.

**ii)** If the function $g(\nabla f, V_i)$, $(i = 2, 3)$ is a non-zero constant along $\alpha$, then $\alpha$ is called a null $f$-eikonal $V_i$-slant helix curve. And, $\nabla f$ is called the axis of the null $f$-eikonal slant helix curve $\alpha$.

**iii)** If the function $g(\nabla f, W)$ is a non-zero constant along $\alpha$, then $\alpha$ is called a null $f$-eikonal Darboux helix curve. And, $\nabla f$ is called the axis of the null $f$-eikonal Darboux helix curve $\alpha$.

**Example 3.1.** We consider the pseudo-Riemannnian manifold $M^3 = \mathbb{R}_1^3$ with the Lorentzian metric $g$. Let consider the function

$$f : M^3 \to \mathbb{R}$$
$$(x, y, z) \to f(x, y, z) = x^2 + y^2 + z$$

given in Example 2.1 and consider the null curve

$$\alpha : I \subset \mathbb{R} \to M^3$$
$$s \to \alpha(s) = (\sinh s, \cosh s, s)$$

in $M^3$. If we compute $\nabla f$, we find out $\nabla f$ as $\nabla f = (2x, 2y, 1)$. Then, we have

$$\|\nabla f\| = \sqrt{\left|1 + 4(-x^2 + y^2)\right|},$$

and, along the curve $\alpha$, we find out

$$\|\nabla f\| = \sqrt{5} = \text{constant}.$$

That is, $f$ is an eikonal function along $\alpha$. Moreover, by a simple computation we have that the tangent and binormal of the curve are

$$V_1(s) = (\cosh s, \sinh s, 1),$$

$$V_2(s) = \left(-\frac{1}{2}\cosh s, \ -\frac{1}{2}\sinh s, \ \frac{1}{2}\right),$$

respectively. Since

$$\nabla f = (2\sinh s, \ 2\cosh s, \ 1),$$

along $\alpha$, we easily see that $g(\nabla f, V_2) = \frac{1}{2} = \text{constant}$ and $g(\nabla f, V_1) = 1 = \text{constant}$ which mean that $\alpha$ is both a null $f$-eikonal helix and $f$-eikonal $V_2$-slant helix in $M^3$.

On the other hand, the curvatures of curve are $\kappa = 1$, $\tau = -\frac{1}{2}$, respectively. Then the Darboux vector of $\alpha$ is

$$W = (0, 0, -1).$$

Then we obtain that

$$g(\nabla f, W) = -1 = \text{constant},$$

along $\alpha$. So, $\alpha$ is also a null $f$-eikonal Darboux helix curve in $M^3$.

Then, we can give the following characterizations for a null curve.

***Theorem 3.1.*** *Let $\alpha : I \subset \mathbb{R} \to M^3$ be a null curve in $M^3$ with non-zero curvatures $\kappa, \tau$ and let $f : M^3 \to \mathbb{R}$ be an eikonal function along curve $\alpha$ and the Hessian $H^f = 0$. If $\alpha(s)$ is a null $f$-eikonal helix curve in $M^3$, then the followings hold,*

*i)* The function $\dfrac{\kappa}{\tau}$ is constant.

*ii)* The axis of null $f$-eikonal helix curve is $\nabla f = c\left(\dfrac{-\tau}{\kappa}V_1 + V_2\right)$, where $g(\nabla f, V_1) = c$ is a non-zero constant.

**Proof.** Let $\alpha(s)$ be a null $f$-eikonal helix curve in $M^3$ with axis $\nabla f$. Then, there exist smooth functions $a_1 = a_1(s)$, $a_2 = a_2(s)$ and $a_3 = a_3(s)$ of arc length $s$ such that

$$\nabla f = a_1 V_1 + a_2 V_2 + a_3 V_3, \tag{22}$$

where $\{V_1, V_2, V_3\}$ is a basis of $TM^3$ (tangent bundle of $M^3$). From (22) we have

$$g(\nabla f, V_1) = a_2 = c = \text{constant}, \quad g(\nabla f, V_2) = a_1, \quad g(\nabla f, V_3) = a_3. \tag{23}$$

Differentiating equalities given in (23), we have $a_3 = 0$, $a_1 = \text{constant}$ and $\dfrac{\kappa}{\tau} = -\dfrac{a_2}{a_1} = \text{constant}$, respectively.

Moreover, from (23) the axis of the null helix is obtained as $\nabla f = c\left(\dfrac{-\tau}{\kappa}V_1 + V_2\right)$, where $g(\nabla f, V_1) = c$ is a non-zero constant.

**Theorem 3.2.** *Let $\alpha: I \subset \mathbb{R} \to \mathbb{R}_1^3$ be a null curve in $\mathbb{R}_1^3$ with non-zero curvatures $\kappa$, $\tau$ and let $f: \mathbb{R}_1^3 \to \mathbb{R}$ be an eikonal function along curve $\alpha$ and the Hessian $H^f = 0$. If $\alpha(s)$ is a null $f$-eikonal $V_2$-slant helix curve in $\mathbb{R}_1^3$, then $\alpha(s)$ is also a null $f$-eikonal helix curve in $\mathbb{R}_1^3$ with axis $\nabla f = c\left(\dfrac{-\tau}{\kappa}V_1 + V_2\right)$ where $g(\nabla f, V_1) = c$ is a non-zero constant.*

**Proof:** Let $\alpha(s)$ be a null $f$-eikonal $V_2$-slant helix curve in $\mathbb{R}_1^3$ with axis $\nabla f$. Then we have $g(\nabla f, V_2) = $ non-zero constant. By differentiation of last equality we get

$$g(\nabla f, V_3) = 0. \tag{24}$$

On the other hand differentiation of $g(\nabla f, V_1)$ in the direction $V_1$ is

$$\nabla_{V_1}[g(\nabla f, V_1)] = \kappa g(\nabla f, V_3),$$

and from (24) we have $g(\nabla f, V_1)$ is a constant. Then $\alpha(s)$ is a null $f$-eikonal helix curve in $\mathbb{R}_1^3$ and from Theorem 3.1, the axis is $\nabla f = c\left(\dfrac{-\tau}{\kappa}V_1 + V_2\right)$ where $g(\nabla f, V_1) = c$ is a non-zero constant.

**Theorem 3.3.** *Let* $\alpha : I \subset \mathbb{R} \to \mathbb{R}_1^3$ *be a null curve in* $\mathbb{R}_1^3$ *with non-zero curvatures* $\kappa, \tau$ *and let* $f : \mathbb{R}_1^3 \to \mathbb{R}$ *be an eikonal function along curve* $\alpha$ *and the Hessian* $H^f = 0$. *If* $\alpha(s)$ *is a null* $f$ *-eikonal helix or* $V_2$ *-slant helix in* $\mathbb{R}_1^3$, *then* $\det(\nabla_{V_1} V_2, \nabla_{V_1}^2 V_2, \nabla_{V_1}^3 V_2) = 0$ *holds.*

**Proof:** Let $\alpha(s)$ be a null curve. Then from Frenet formulae (21) we have the followings

$$\begin{cases} \nabla_{V_1}^2 V_2 = -\tau^2 V_1 - \tau \kappa V_2 + \tau' V_3, \\ \nabla_{V_1}^3 V_2 = -3\tau\tau' V_1 - \left((\tau\kappa)' + \kappa\tau'\right) V_2 + (-2\kappa\tau^2 + \tau'') V_3. \end{cases} \qquad (25)$$

From (25) we have $\det(\nabla_{V_1} V_2, \nabla_{V_1}^2 V_2, \nabla_{V_1}^3 V_2) = \tau^5 \left(\dfrac{\kappa}{\tau}\right)'$. Then by Theorem 3.2 and Theorem 3.1, we say that if $\alpha(s)$ is a null $f$-eikonal helix or $V_2$-slant helix curve in $\mathbb{R}_1^3$ then $\det(\nabla_{V_1} V_2, \nabla_{V_1}^2 V_2, \nabla_{V_1}^3 V_2) = 0$ holds.

**Theorem 3.4.** *Let* $\alpha : I \subset \mathbb{R} \to M^3$ *be a null curve in* $M^3$ *with curvatures* $\kappa, \tau$ *and let* $f : M^3 \to \mathbb{R}$ *be an eikonal function along curve* $\alpha$ *and the Hessian* $H^f = 0$. *If* $\alpha(s)$ *is a null* $f$ *-eikonal* $V_3$ *-slant helix curve in* $M^3$, *then the following properties hold:*

*i)* $\kappa(s) \displaystyle\int_0^s \tau(s) ds + \tau(s) \displaystyle\int_0^s \kappa(s) ds = 0$ *holds, where* $g(\nabla f, V_3) = c$ *is a non-zero constant.*

*ii) The axis of the* $V_3$ *-slant helix is given by*

$$\nabla f = c \left[ \left( \int_0^s \tau(s) \, ds \right) V_1 + \left( \int_0^s \kappa(s) \, ds \right) V_2 + V_3 \right] \qquad (26)$$

**Proof:** Since we assume that $\alpha(s)$ is a null $f$-eikonal $V_3$-slant helix curve in $M^3$, we have $g(\nabla f, V_3) = c$ is a non-zero constant. Then we can write

$$\nabla f = a_1(s) V_1 + a_2(s) V_2 + c V_3, \qquad (27)$$

where $a_i = a_i(s); (i = 1, 2)$ are the differentiable functions of $s$. From Corollary 1.1, $\nabla f$ is parallel in $M^3$, i.e., $\nabla_{V_1} \nabla f = 0$ along $\alpha$. Then from (27) we obtain

$$(a_1' - c\tau) V_1 + (a_2' - c\kappa) V_2 + (a_1 \kappa + a_2 \tau) V_3 = 0, \qquad (28)$$

which gives the following system

$$a_1' - c\tau = 0, \quad a_2' - c\kappa = 0, \quad a_1 \kappa + a_2 \tau = 0. \qquad (29)$$

And from (29) and (27), we have the followings immediately,

$$\nabla f = c\left[\left(\int_0^s \tau(s)\,ds\right)V_1 + \left(\int_0^s \kappa(s)\,ds\right)V_2 + V_3\right],$$

$$\kappa(s)\int_0^s \tau(s)ds + \tau(s)\int_0^s \kappa(s)ds = 0.$$

Theorem 3.4 gives us the following corollary:

***Corollary 3.1.*** *Let $\alpha : I \subset \mathbb{R} \to M^3$ be a null curve in $M^3$ with curvatures $\kappa, \tau$ and let $f : M^3 \to \mathbb{R}$ be an eikonal function along curve $\alpha$ and the Hessian $H^f = 0$. If $\alpha(s)$ is a null $f$-eikonal $V_3$-slant helix curve in $M^3$. Then the followings holds,*

i) *$\alpha(s)$ is a null $f$-eikonal helix curve if and only if $\kappa(s) = 0$.*

ii) *$\alpha(s)$ is a null $f$-eikonal $V_2$-slant helix curve if and only if $\tau(s) = 0$.*

**Proof:** From Theorem 3.4, we have that the axis is given by

$$\nabla f = a_1(s)V_1 + a_2(s)V_2 + cV_3,$$

where $c$ is a non-zero constant. Then we have that

$$g(\nabla f, V_1) = a_2(s),\ \ g(\nabla f, V_2) = a_1(s), \tag{30}$$

and from (29) and (30) we have the followings,

i) $g(\nabla f, V_1) = a_2$ is constant if and only if $\kappa(s) = 0$,

ii) $g(\nabla f, V_2) = a_1$ is constant if and only if $\tau(s) = 0$,

which finish the proof.

***Theorem 3.5.*** *Let $\alpha : I \subset \mathbb{R} \to M^3$ be a null $f$-eikonal Darboux helix with Darboux vector $W = \tau V_1 - \kappa V_2$ in $M^3$ with non-zero curvatures $\kappa, \tau$, where $f : M^3 \to \mathbb{R}$ is an eikonal function along curve $\alpha$ and the Hessian $H^f = 0$. Then $\alpha$ is a null $f$-eikonal $V_3$-slant helix if and only if $\kappa\tau$ is constant.*

**Proof.** Since $\alpha$ is a null $f$-eikonal Darboux helix, we have $g(W, \nabla f) = $ constant. On the other hand, there exist smooth functions $a_1 = a_1(s)$, $a_2 = a_2(s)$ and $a_3 = a_3(s)$ of arc length $s$ such that

$$\nabla f = a_1 V_1 + a_2 V_2 + a_3 V_3, \tag{31}$$

where $a_1, a_2, a_3$ are assumed non-zero and $\{V_1, V_2, V_3\}$ is a basis of $TM^3$. Since $\nabla_{V_1} \nabla f = 0$ along $\alpha$, if we take the derivative in each part of (31) in the direction $V_1$ in $M^3$ and use the Frenet equations, we get

$$(a_1' - a_3\tau)V_1 + (a_2' - a_3\kappa)V_2 + (a_3' + a_1\kappa + a_2\tau)V_3 = 0, \tag{32}$$

where $a_i'(s) = V_1[a_i]$, $(i = 1, 2, 3)$. Since the Frenet frame $\{V_1, V_2, V_3\}$ is linearly independent, we have the system

$$\begin{cases} a_1' - a_3\tau = 0, \\ a_2' - a_3\kappa = 0, \\ a_3' + a_1\kappa + a_2\tau = 0. \end{cases} \tag{33}$$

Equality $g(W, \nabla f) = \text{constant}$ gives that

$$a_2\tau - a_1\kappa = \text{constant}. \tag{34}$$

Differentiating (34) and using the first and second equations of system (33) we obtain

$$a_2\tau' - a_1\kappa' = 0. \tag{35}$$

From (35) and the third equation of system (33) it follows

$$a_3' = -\frac{(\kappa\tau)'}{\kappa'} a_2. \tag{36}$$

If $a_2 = 0$ in (36), from (33) we see that $a_1 = a_3 = 0$ which is a contradiction. Then $a_2 \neq 0$ and we have that $a_3 = a_3(s)$ is a constant if and only if $\kappa\tau$ is constant, which means that $\alpha$ is a null $f$-eikonal $V_3$-slant helix if and only if $\kappa\tau$ is constant.

## References

[1] Ali, A.T., Position vectors of slant helices in Euclidean Space $E^3$, J. of Egyptian Math. Soc., 20(1) (2012) 1-6.
[2] Ali, A.T., Lopez, R., Slant Helices in Minkowski Space $E_1^3$, J. Korean Math. Soc. 48(1) (2011) 159-167.
[3] Ali, A.T., Turgut, M., Position vector of a time-like slant helix in Minkowski 3-space, J. Math. Anal. Appl. 365 (2010) 559–569.
[4] Barros, M., General helices and a theorem of Lancret, Proc. Amer. Math. Soc. 125, no.5, (1997) 1503–1509.
[5] Di Scala, A.J., Ruiz-Hernandez, G., Higher codimensional euclidean helix submanifolds, Kodai Math. J. 33, (2010) 192-210.
[6] Duggal, K.L., Jin, D.H., Null curves and hypersurfaces of semi-Riemannian manifolds, World Scientific, 2007.
[7] Ekmekçi, N., Hacısalihoğlu, H.H., On Helices of a Lorentzian Manifold, Commun. Fac. Sci. Univ. Ank. Series A1, 45 (1996) 45-50.
[8] Ferrandez, A., Gimenez, A., Lucas, P., Null generalized helices in Lorentz-Minkowski spaces, J. Phys. A 35, no.39 (2002) 8243-8251.

[9] Ferrandez, A., Gimenez, A., Lucas, P., Null helices in Lorentzian space forms, Internat. J. Modern Phys. A. 16, no.30 (2001) 4845-4863.

[10] Fischer, A.E., Riemannian maps between Riemannian manifolds, Contemporary Math., Vol 182, (1992) 331-366.

[11] Gök, İ., Camcı, Ç, Hacısalihoğlu, H.H., $V_n$-slant helices in Euclidean $n-$space $E^n$, Math. Commun., 14(2) (2009) 317-329.

[12] Izumiya, S., Takıyama, A., A time-like surface in Minkowski 3-space which contain pseudo-circles, Proc. Edinburg Math. Soc., 40 (1997), 127-136.

[13] Izumiya, S., Takiyama, A., A time-like surface in Minkowski 3-space which contain light-like lines, Journal of Geometry, 64 (1999), 95-101.

[14] Izumiya, S., Takeuchi, N., New special curves and developable surfaces, Turk J. Math. 28, (2004) 153-163.

[15] Kobayashi, O., Maximal surfaces in the 3-dimensional Minkowski space $L^3$, Tokyo J. Math., 6 (1983), 297-309.

[16] Kula, L. and Yayli, Y., On slant helix and its spherical indicatrix, Appl. Math. and Comp., 169, (2005) 600-607.

[17] Kula, L., Ekmekçi, N., Yaylı, Y., İlarslan, K., Characterizations of Slant Helices in Euclidean 3-Space, Turk J Math., 33 (2009) 1–13.

[18] O'Neill, B., Semi-Riemannian Geometry with Applications to Relativity. Academic Press, London (1983).

[19] Önder, M., Kazaz, M., Kocayiğit, H. Kılıç, O., $B_2$-slant helix in Euclidean 4-space $E^4$, Int. J. Contemp. Math. Sci. 3(29-32) (2008) 1433-1440.

[20] Struik, D.J., Lectures on Classical Differential Geometry, 2nd ed. Addison Wesley, Dover, (1988).

[21] Şenol, A., Zıplar, E., Yaylı, Y., On $f$-Eikonal Helices and $f$-Eikonal Slant Helices in Riemannian Manifolds, arXiv:1211.4960 [math.DG].

[22] Zıplar, E., Şenol, A., Yaylı, Y., On Darboux helices in Euclidean 3-space, Global J. of Sci. Frontier Res., 12(13) (2012) 72-80.